\documentclass[11pt]{amsart}
\usepackage[english]{babel}
\usepackage{amssymb,amscd}

\usepackage[all]{xy}

\usepackage{euscript}
\usepackage{amsmath}
\usepackage{ifpdf}

\usepackage{amsfonts}
\usepackage{mathrsfs}

\usepackage{amsmath}

\newtheorem{Theorem}{Theorem}[section]
\newtheorem{Lemma}[Theorem]{Lemma}
\newtheorem{Corollary}[Theorem]{Corollary}
\newtheorem{Proposition}[Theorem]{Proposition}

\newtheorem{Remark}[Theorem]{Remark}
\newtheorem{Example}[Theorem]{Example}

\newtheorem{Definition}[Theorem]{Definition}
\newtheorem{notation}[Theorem]{Notation}
\newtheorem{Question}[Theorem]{Question}

\def\sqr#1#2{{\vcenter{\hrule height.#2pt
\hbox{\vrule width.#2pt height#1pt \kern#1pt \vrule width.#2pt}
\hrule height.#2pt}}}
\def\proof{\noindent{\bf Proof. }}
\def\qed{\hspace*{\fill} $\square$}

\def\fa{{\mathfrak a}}

\def\fm{{\mathfrak m}}

\def\fp{{\mathfrak p}}

\def\Cc{\mathscr{C}}

\def\Hc{\mathscr{H}}

\def\Jc{\mathscr{J}}

\def\Kc{\mathscr{K}}

\DeclareMathOperator{\Ann}{ann}
\DeclareMathOperator{\Hom}{Hom}

\DeclareMathOperator{\height}{height}

\begin{document}

\title{ Strong $F$-regularity and generating morphisms of local cohomology modules}

\author{Mordechai Katzman}
\address{Department of Pure Mathematics,
University of Sheffield, Hicks Building, Sheffield S3 7RH, United Kingdom}
\email{M.Katzman@sheffield.ac.uk}

\author{Cleto B.~Miranda-Neto}
\address{Departamento de Matem\'atica, Universidade Federal da
Para\'iba, 58051-900 Jo\~ao Pessoa, PB, Brazil.}
\email{cleto@mat.ufpb.br}

\thanks{
C. B. Miranda-Neto thanks the Department of Pure Mathematics, University of Sheffield (United Kingdom), for hospitality during his 1-year stay as a visiting researcher. He was partially supported by CAPES-Brazil 
(``Est\'agio S\^enior no Exterior", grant 88881.121012/2016-01), and by CNPq-Brazil (``Chamada Universal 2016", grant 421440/2016-3)}


\begin{abstract}

  We establish a criterion for the strong $F$-regularity of a (non-Gorenstein) Cohen-Macaulay reduced complete local ring of dimension at least $2$, 
  containing a perfect field of prime characteristic $p$.
  We also describe an explicit generating morphism (in the sense of Lyubeznik) for the top local cohomology module with support in certain ideals arising from an $n\times (n-1)$ matrix $X$ of indeterminates. For $p\geq 5$, these results led us to derive a simple, new proof of the well-known fact that the generic determinantal ring defined by the maximal minors of $X$ is strongly $F$-regular.

\end{abstract}

\maketitle

\section{Introduction}

In this paper we deal with the property of strong $F$-regularity as well as the detection of an explicit generating morphism (in the sense of Lyubeznik \cite{L}) for certain local cohomology modules. 
As the main application of our results, we recover, in a quite simple way, the fact that the determinantal ring $A$ defined by the maximal minors of an $n\times (n-1)$ generic matrix over a perfect field of 
prime characteristic $p\geq 5$
is strongly $F$-regular, a result that was proved, in more generality, first by Hochster and Huneke \cite{HH} and later by Bruns and Conca \cite{BC} via completely different methods. 

It is worth mentioning that Conca and Herzog \cite{CH}, also Glassbrenner and Smith \cite{G-S}, proved that certain types of ladder determinantal rings are strongly $F$-regular (hence $F$-rational and $F$-pure as well). Moreover, the $F$-purity of Hankel determinantal rings has been verified in \cite{CMSV}.

The central result of the present paper, Theorem \ref{main}, provides an effective test for the strong $F$-regularity of a (non-Gorenstein) Cohen-Macaulay reduced complete local domain $A$ of dimension at least $2$, containing a perfect field of positive characteristic.
As we point out in Remark \ref{pure}, it affords a structural similarity with the characterization given in \cite{Glass}. 
We apply the theorem in the case where $A$ is the determinantal ring defined by the minors of order $2$ of a $3\times 3$ generic symmetric matrix over a field of characteristic $2$ or $3$ (Example \ref{symm}). In Example \ref{N2} we illustrate that our result may fail if we relax the main hypothesis. Furthermore, we show in Example \ref{false-converse} that the converse of our theorem is not true.

Another main result (Proposition \ref{prod}) gives an explicit generating morphism for the local cohomology module
$$H^2_{I_{n-1}(X)}(k[[X]])$$ where $X$ is an $n\times (n-1)$ matrix of indeterminates over a field $k$ of positive characteristic, for any $n\geq 2$.
As it turned out, we prove in Remark \ref{u'} that an adaptation of the proof also yields a generating morphism for the module
$$H^3_{I_{n-1}(X)+I_{n-2}(X^\prime)}(k[[X]])$$ whenever $n\geq 3$, where $X^\prime$ is the submatrix formed by the first $n-2$ rows of $X$.
Needless to say, the study of local cohomology supported at determinantal ideals has been a major topic of research under various viewpoints
(see for example  \cite{LSW},  \cite{RW},  \cite{RWW} and  \cite{W}). Our result is interesting in its own in view of Lyubeznik's theory \cite{L},
and moreover it will serve as a crucial ingredient in the main application furnished in this paper, regarding the strong $F$-regularity of certain rings.
More precisely, from Theorem \ref{main} and Remark \ref{u'} we derive, in Corollary \ref{generic}, a new and quite simple proof of the well-known fact mentioned in the first paragraph: the generic determinantal ring $k[[X]]/I_{n-1}(X)$ is strongly $F$-regular if the ground field $k$ is perfect and of characteristic $\geq 5$.

We close the paper with a few observations and questions on the extension of our methods to the case of arbitrary (not necessarily maximal) minors
as well as to other types of determinantal rings, such as those defined by minors of a generic \emph{symmetric} matrix.
While the structure result given in Proposition \ref{prod} and the related fact observed in Remark \ref{u'} are no longer valid in these situations,
it would be interesting to explore adapting our methods with a view to tackling these new cases.

\section{Preliminaries: notions and methods in prime characteristic}\label{Section: Preliminaries}

In this section we provide a brief review on the main concepts and methods in positive characteristic that we will use in this paper.
Unless explicitly stated otherwise, by the term {\it ring} we shall mean, in the entire paper, {\it Noetherian commutative unital ring}.

Let $A$ be a ring of prime characteristic $p$. For an integer $e\geq 0$, let $f^e\colon A\rightarrow A$ be the $e$th Frobenius endomorphism of $A$, i.e., the function $f^e(a)=a^q$ for $a\in A$ and $q=p^e$.

Given an ideal $\fa \subset A$, we denote by $\fa^{[q]}$ the ideal generated by the $q$th powers of all elements of $\fa$.
Clearly, if $\fa=(x_1,\ldots, x_n)$ then $$\fa^{[q]}~=~(x_1^q,\, \ldots, \, x_n^q).$$
\begin{Definition}\rm
The {\it tight closure} of an ideal $\fa\subseteq A$ is the ideal $\fa^*$ consisting of $x\in A$ such that $cx^q\in \fa^{[q]}$ for some $c\in A$ not in a minimal prime of $A$, and all large $q=p^e$.
\end{Definition}
It turns out that $\fa^*$ is an ideal satisfying $\fa \subseteq \fa^*=(\fa^*)^*$, and $\fa$ is said to be {\it tightly closed} if $\fa=\fa^*$.
For details on this closure operation we refer to \cite{HH-1} and \cite{Hu}.

Given an $A$-module $M$, we can give it  a new $A$-module structure via $f^e$. To this end,
let $F^e_* M$ stand for the additive Abelian group $M$, with elements denoted by
$F^e_* a$,\, $a\in M$, and endow $F^e_* M$  with the $A$-module structure given by  $a F_*^e m = F^e_* a^{p^e} m$.

\begin{Definition}\rm (cf.~\cite{HH} and \cite[Section 5]{HH0})
\label{Definition: Fregularity}
The ring $A$ is said to be
\begin{enumerate}
 \item[(a)] \emph{strongly $F$-regular}, if for each non-zero $c\in A$, the $A$-linear map $A \rightarrow F_*^e A$ sending
 $1$ to $F_*^e c$ splits for all large $e$;
 \item[(b)] \emph{weakly $F$-regular}, if every ideal in $A$ is tightly closed;
 \item[(c)] \emph{$F$-rational}, if every parameter ideal in every localization of $A$ is tightly closed;
 \item[(d)] \emph{$F$-pure}, if for all $A$-modules $M$, the map $f\otimes_A 1 : A\otimes_A M \rightarrow F_*^1 A \otimes M$ is injective.
\end{enumerate}
\end{Definition}
 It is well-known that (a) $\Rightarrow$ (b) $\Rightarrow$ (c), (b) $\Rightarrow$ (d) (cf.~\cite{F-W}), and (c) implies that $A$ is a Cohen-Macaulay normal domain (cf.~\cite{HH0}).

Given any $A$-linear map $g: M \rightarrow F_*^e M$, we have an additive map  $\widetilde{g} : M \rightarrow M$ obtained by identifying $F_*^e M$ with $M$.
We point out that $\widetilde{g}$ is not $A$-linear; instead, it satisfies $\widetilde{g}(a m)=a^{p^e} \widetilde{g}(m)$ for all $a\in A$ and $m\in M$.
We call additive maps with this property \emph{$e$th Frobenius maps}.
Conversely, a Frobenius map $h: M \rightarrow M$ defines an $A$-linear map  $M \rightarrow F_*^e M$ given by $m\mapsto F_*^e h(m)$.

To keep track of Frobenius maps we introduce the following  skew-polynomial rings.
Let $A[\Theta; f^e]$ be the free $A$-module $\bigoplus_{i\geq 0} A \Theta^i$ and
give $A[\Theta; f^e]$ the structure of a ring by defining the (non-commutative) product $(a \Theta^i) (b \Theta^j) = a b^{p^{e i}} \Theta^{i+j}$
(see~\cite[Chapter 1]{LamAFirstCourseInNoncommutativeRings} for a more general version of this construction).
Now an $e$th Frobenius map $\widetilde{g}$  on an $A$-module $M$ corresponds to an
$A[\Theta; f^e]$-module structure on $M$ where the action of $\Theta$ on $M$ is given by $\Theta m = \widetilde{g}(m)$ for all $m\in M$.

A crucial set of tools in the prime characteristic toolkit are the \emph{Frobenius functors} which we now define.
For any $A$-module $M$ and $e\geq 1$, we can extend scalars and obtain the $F_*^e A$-module  $F_*^e A \otimes_A M$. If we now identify the rings $A$ and
$F_*^e A$, we obtain the $A$-module $A \otimes_R M$ where
for $a, b\in A$ and $m\in M$, $a(b\otimes m)=ab\otimes m$ and $a^{p^e} b \otimes m=b\otimes a m$, and we denote this module by $F^e_A(M)$.
Clearly, homomorphisms $M\rightarrow N$ induce $A$-linear maps $F_A^e(M) \rightarrow F_A^e(M)$ and thus we obtain the
\emph{$e$th Frobenius functors} $F_A^e(-)$. When the ring $A$ is regular, a classical result due to Kunz (cf.~\cite{Kunz}), which we shall use tacitly in this paper, states that the functor $F^e_A(-)$ is exact. A useful consequence is that in this case, if $\fa \subset A$ is an ideal, then the $A$-modules $A/\fa$ and $A/\fa^{[p]}$ possess the same set of associated primes (cf.~\cite[Proposition 21.11]{24HoursLocalCohomology}).

An $e$th Frobenius map $g: M \rightarrow M$
gives rise to an $A$-linear map $ \overline{g} : F_R^e (M) \rightarrow M$ defined by $\overline{g}(r\otimes m)=r g(m)$;
this is well defined since for all $a, b \in A$ and $m\in M$
$$\overline{g}(a^{p^e} b\otimes m) = a^{p^e} b g(m) =  b g(a m) =  \overline{g}(b \otimes am) .$$
In the special case when $M$ is an Artinian module over a complete regular ring, this gives a way to define a
``Matlis-dual which keeps track of Frobenius'' functor, which we describe next.

Henceforth in this paper we adopt the following notation.

\begin{notation}\label{Notation: complete regular local}\rm
Let $(R, \mathfrak{m})$  be a $d$-dimensional complete regular local ring of prime characteristic $p$ and let $A$
be its quotient by an ideal $I\subset R$.
We will denote $E=E_R(R/\mathfrak{m})$ and $E_A=E_A(A/\mathfrak{m}A)=\Ann_E I$
the injective hulls of the residue fields of $R$ and $A$, respectively.
The Matlis dual functor $\Hom_R (-, E)$ will be denoted $(-)^\vee$.
\end{notation}

A crucial ingredient for the construction that follows is the fact that for both Artinian and Noetherian modules $R$-modules $M$,
there is a natural identification of $F_R^e(M)^\vee$ with $F_R^e(M^\vee)$ (\cite[Lemma 4.1]{L})
and  henceforth we identify these tacitly.
A map of Artinian $R[\Theta; f^e]$-modules $\rho: M\rightarrow N$
yields a commutative diagram
\begin{equation*}
\xymatrix{
F_R^e(M) \ar@{>}[r]^{F_R^e (\rho)} \ar@{>}[d]^{1\otimes \Theta} & F_R^e(N) \ar@{>}[d]^{1\otimes \Theta} \\
M \ar@{>}[r]^{\rho} & N \\
}
\end{equation*}
and an application of the Matlis dual gives the commutative diagram
\begin{equation*}
\xymatrix{
N^\vee \ar@{>}[r]^{\rho^\vee} \ar@{>}[d]^{1\otimes \Theta^\vee} & M^\vee \ar@{>}[d]^{1\otimes \Theta^\vee}  \\
F_R^e(N)^\vee \ar@{>}[r]^{F_R^e (\rho^\vee)} & F_R^e(M)^\vee \\
}
\end{equation*}
Define $\mathcal{C}_e$ to be the category of
Artinian $R[\Theta; f^e]$-modules and $\mathcal{D}_e$ the category of $R$-linear maps $N \rightarrow F^e_R (N)$ for Noetherian $R$-modules $N$,
where morphisms in $\mathcal{D}_e$ are commutative diagrams
\begin{equation}\label{CD1}
\xymatrix{
N \ar@{>}[r]^{\varphi} \ar@{>}[d]^{\xi} & M \ar@{>}[d]^{\zeta}  \\
F_R^e(N) \ar@{>}[r]^{F_R^e(\varphi)} & F_R^e(M) \\
}
\end{equation}
The construction above yields a contravariant functor $\Delta^e : \mathcal{C}_e \rightarrow \mathcal{D}_e$, and this functor is exact \cite[Chapter 10]{BS}.
Furthermore, an application of the Matlis dual to (\ref{CD1}) yields
\begin{equation}\label{CD2}
\xymatrix{
F_R^e(M^\vee) \ar@{>}[r]^{F_R^e(\varphi^\vee)} \ar@{>}[d]^{\zeta^\vee} & F_R^e(N^\vee) \ar@{>}[d]^{\xi^\vee}\\
M^\vee \ar@{>}[r]^{\varphi^\vee}  & N^\vee   \\
}.
\end{equation}
which can be used to equip $M^\vee$ and $N^\vee$ with $R[\Theta; f^e]$-module structures given by
$\Theta m=\zeta^\vee (1\otimes m)$ and $\Theta n=\xi^\vee (1\otimes n)$, respectively.
With these structures, $\varphi^\vee$ is $R[\Theta; f^e]$-linear.
This construction yields an exact contravariant functor $\Psi^e : \mathcal{D}_e \rightarrow \mathcal{C}_e$. Now, after the identification of the double Matlis dual $(-)^{\vee\vee}$ with the identity
functor on Artinian and Noetherian $R$-modules, the compositions $\Psi^e \circ \Delta^e$ and $\Delta^e \circ \Psi^e$ yield the identity functors on
$\mathcal{C}_e$ and $\mathcal{D}_e$, respectively (cf.~\cite{K} for details).
In this sense, we can think of $\Delta^e$ as the Matlis dual that keeps track of a given Frobenius map.

In this paper we will focus our attention on a specific family of Artinian modules with Frobenius maps, namely, the local cohomology modules
$H^{\dim R}_\mathfrak{m}(R) \cong E$ and $H^\bullet_\mathfrak{m}(A)$, their submodules and homomorphic images.
Recall that for any ideal $J$ in a commutative ring $S$ of prime characteristic, the Frobenius map $f:S \rightarrow S$ induces a Frobenius map
$H^i_J(S) \rightarrow H^i_J(S)$ (\cite[Chapter 21]{24HoursLocalCohomology}), and when $J$ is a maximal ideal, these local cohomology modules are Artinian.

Finally, we recall some of the constructions introduced in Lyubeznik \cite{L}.
An $R$-module $\mathcal{M}$ together with an isomorphism $\theta: \mathcal{M} \rightarrow F_R(\mathcal{M})$ is an \emph{$F$-module with structural isomorphism $\theta$}.
If $M$ is an $R$-module and $\phi :  M \rightarrow F_R(M)$ is $R$-linear, we obtain the following $F$-module as a direct limit
$$\mathcal{M} \, = \, \displaystyle \lim_{\longrightarrow} ~ (M \xrightarrow[]{\phi} F_R(M) \xrightarrow[]{F_R(\phi)} F^2_R(M) \xrightarrow[]{F^2_R(\phi)} \, \cdots)$$
We call such an $F$-module a \emph{$F$-finite $F$-module with generating morphism $\phi$} (if $\phi$ is injective, we call $\phi$ a {\it root} of $\mathcal{M}$).

If $M$ is Artinian as an $R[\Theta; f]$-module,
there is a natural $R$-linear map $\alpha_M \colon F_R(M)\rightarrow M$, given by $\alpha_m( r \otimes m)= r \Theta m$
with Matlis dual map $\alpha_M^{\vee}\colon M^{\vee}\rightarrow F_R(M)^{\vee}\cong F_R(M^{\vee})$.
We can now define
{\it Lyubeznik functor} as
$$\Hc_{R, A}(M) \, = \, \displaystyle \lim_{\longrightarrow} ~ (M^{\vee} \xrightarrow[]{\alpha_M} F_R(M^{\vee})\xrightarrow[]{F_R(\alpha_M)} F_R^2(M^{\vee})\xrightarrow[]{F_R^2(\alpha_M)} \, \cdots)$$
which is an exact, contravariant functor from the category of Artinian $R[\Theta; f]$-modules to the category of $F$-finite $F$-modules.

\section{Generating morphisms of local cohomology modules}\label{gen}

In this section we will further assume that the ring $A=R/I$ is reduced, Cohen-Macaulay of dimension $d\geq 1$, and will write $h=\height I$. Our objective is to describe some generating morphisms for $H_{I}^h(R)$, the only non-vanishing local cohomology module with support at $I$ (cf.~\cite[Proposition III.4.1]{PS}).

First we recall a basic fact:

\begin{Lemma}{\rm (\cite[Proposition 3.3.18]{BH})}
\label{can}
Let $B$ be a Cohen-Macaulay ring having a canonical module $\omega_B$.
If $B$ is generically Gorenstein then $\omega_B$ can be identified with an ideal $K\subset B$. For any such identification, either $K$ has height $1$ $($in which case $B/K$ is Gorenstein$)$ or $K = B$.
\end{Lemma}

Since our ring $A$ is Cohen-Macaulay and reduced (hence, generically Gorenstein), Lemma \ref{can} implies that the canonical module $\omega = \omega_A$ of $A$ can be identified with an ideal
$\Omega/I\subset A$, for a suitable ideal $\Omega \supset I$ which has height $h+1$ if $\Omega \neq R$.

We also invoke the following observation:

\begin{Lemma}{\rm (\cite[Example 3.7]{LS})}\label{cyclic}
If $(B, {\mathfrak n})$ is a complete Cohen-Macaulay local ring of dimension $b\geq 1$ containing a field of positive characteristic,
then the $B$-module of all Frobenius maps on $H^{b}_{{\mathfrak n}}(B)$ is free of rank $1$, generated by the natural map induced by the Frobenius endomorphism of $B$.
\end{Lemma}

Next we furnish a result that will be extremely useful in the sequel, in particular for our Theorem \ref{main} (where we furnish a criterion of strong $F$-regularity). Further details, even in more generality, 
can be found in  \cite[Subsection 3.4.2]{Bli},  \cite[Subsection 5.3]{KMSZ} and \cite[Section 4]{L}, but we supply a proof herein for the reader's convenience.

\begin{Proposition}\label{useful} Assume that $d\geq 1$.
\begin{enumerate}

\item[{\rm (i)}] We can identify\, $H^{d}_{{\mathfrak m}}(A)  = \Ann_{E}(I)/\Ann_{E}(\Omega)$;

\item[{\rm (ii)}] The $R$-module of all Frobenius maps on $H^{d}_{{\mathfrak m}}(A)$ is isomorphic to $$((I^{[p]}\colon I)\cap (\Omega^{[p]}\colon \Omega))/I^{[p]}$$

\item[{\rm (iii)}] Under the isomorphism in (ii), the natural Frobenius map on $\Ann_{E}(I)/\Ann_{E}(\Omega)$ is given 
(up to a unit) by $u T$, where $T$ is the natural Frobenius map on $E$, $u\in (I^{[p]}\colon I)\cap (\Omega^{[p]}\colon \Omega)$,
and the image of $u$ in $((I^{[p]}\colon I)\cap (\Omega^{[p]}\colon \Omega))/I^{[p]}$ generates this module.

\end{enumerate}
\end{Proposition}
\proof The inclusion $\omega \subset A$ is compatible with the Frobenius endomorphism $f\colon A\rightarrow A$, and the natural short exact sequence
$$0\longrightarrow \omega \longrightarrow A \longrightarrow A/\omega \longrightarrow 0$$ induces an exact sequence of $A[\Theta; f]$-modules
$$H^{d-1}_{{\mathfrak m}}(A)\longrightarrow H^{d-1}_{{\mathfrak m}}(A/\omega) \longrightarrow  H^{d}_{{\mathfrak m}}(\omega) \longrightarrow  H^{d}_{{\mathfrak m}}(A) \longrightarrow  H^{d}_{{\mathfrak m}}(A/\omega).$$
Since $A$ is Cohen-Macaulay, $H^{d-1}_{{\mathfrak m}}(A)=0$, and also $H^{d}_{{\mathfrak m}}(A/\omega)=0$ as $\dim(A/\omega)<d$
(or trivially if $\omega = A$). Moreover, we can identify $H^{d}_{{\mathfrak m}}(\omega)=E_A=\Ann_{E}(I)$.
Since $A/\omega\cong R/\Omega$ is Gorenstein (cf.~\cite[Proposition 3.3.11(b)]{BH}), we get that
$H^{d-1}_{{\mathfrak m}}(A/\omega)$ is the injective hull of the residue field of $A/\omega$, hence
the annihilator of
$H^{d-1}_{{\mathfrak m}}(A/\omega)$ is $\omega = \Omega/I$, and we may write $H^{d-1}_{{\mathfrak m}}(A/\omega)=\Ann_{E}(\Omega)$. It follows
a short exact sequence
\begin{equation}\label{eq1}
0\longrightarrow  \Ann_{E}(\Omega)  \longrightarrow  \Ann_{E}(I)   \longrightarrow   H^{d}_{{\mathfrak m}}(A)  \longrightarrow 0
\end{equation}
where the injection is an inclusion. This gives (i).

\smallskip

An application of the functor $\Delta^1$ described in Section \ref{Section: Preliminaries} to the short exact sequence of $R[\Theta; f]$-modules (\ref{eq1})
yields a short exact sequence in $\mathcal{D}^1$
\begin{equation}\label{eq2}
\xymatrix{
0 \ar@{>}[r] 	&	\Omega/I \ar@{>}[r] \ar@{>}[d]^{u} 	& 	R/I \ar@{>}[r] \ar@{>}[d]^{u}	& 	R/\Omega \ar@{>}[r] \ar@{>}[d]^{u}	&   0 \\
0 \ar@{>}[r] 	&	\Omega^{[p]}/I^{[p]} \ar@{>}[r] 	& 	R/I^{[p]} \ar@{>}[r]		& 	R/\Omega^{[p]} \ar@{>}[r] 		&   0 \\
}
\end{equation}
where the vertical maps are multiplication by some $u\in R$. In order for these maps to be well-defined, we must have
$u\in (I^{[p]}\colon I)\cap (\Omega^{[p]}\colon \Omega)$. Since $\Omega/I$ must contain a non-zero-divisor, the left-most vertical map is zero if and only if $u\in I^{[p]}$.
Conversely, given any $u\in R$ that makes (\ref{eq2}) commute, an application of the functor $\Psi^1$ will endow the modules in (\ref{eq1}) with an $R[\Theta; f]$-module structure making the maps there
$R[\Theta; f]$-linear.

As for (iii), Lemma \ref{cyclic} shows that the natural Frobenius map on $H^d_{\mathfrak{m}}(A)$ generates the module of all Frobenius maps on  $H^d_{\mathfrak{m}}(A)$. \qed

\begin{Remark}\label{generating}\rm From the proof above and the well-known description of $H_I^h(R)$ as the $F$-module $\Hc_{R, A}(H^d_{\mathfrak{m}}(A))$, it follows that if $u\in R$ is such that its image modulo $I^{[p]}$ generates the cyclic module $((I^{[p]}\colon I)\cap (\Omega^{[p]}\colon \Omega))/I^{[p]}$, then the (well-defined) multiplication map
$$\Omega/I \stackrel{\cdot u}{\longrightarrow} \Omega^{[p]}/I^{[p]}$$
is isomorphic to a generating morphism of $H_I^h(R)$.

\end{Remark}

\section{A general criterion of strong $F$-regularity}\label{F-rat}

As in the previous section, we let $A=R/I$ be a Cohen-Macaulay reduced ring of dimension $d\geq 1$ (if $d=0$ then $A$ is simply a field), where $(R, {\mathfrak m})$ is a formal power
series ring over a field $k$ of characteristic $p>0$, which, moreover, we shall tacitly assume to be a perfect field, and we let $\Omega/I$ be a canonical ideal of $A$.

In accordance with standard terminology, we say that a radical ideal $J\supset I$ {\it defines the singular locus} of $A$ if\, ${\rm Sing}(A)=V(J)$, i.e., for a prime ideal $\fp \supset I$, the local ring $A_{\fp/I}$ is regular if and only if $\fp \nsupseteq J$.  In order to identify such an ideal in general, we may simply resort to the classical method by means of (formal) derivatives, namely, if $\Jc$ stands for the Jacobian matrix of some (in fact, any) set of generators of $I$, then we can take $$J\, =\, \sqrt{I_h(\Jc)\, +\, I}$$ where $h$ is the height of $I$ and $I_h(\Jc)$ is the ideal generated by the subdeterminants of $\Jc$ of order $h$.

We consider, in addition, three more definitions (details on the first two notions can be found, under different terminology, in  \cite{K}).

\begin{Definition}\label{u}\rm Fix some $v\in R$. An ideal $L\subset R$ is said to be {\it $v$-compatible} if\, $v L  \subset  L^{[p]}$.

\end{Definition}

\noindent Note that $L$ is $v$-compatible if and only if  $\Ann_{E}(L)$ is invariant under $v T$, where $T\colon E\rightarrow E$ is the natural Frobenius map. In this case, the map $$R/L \stackrel{\cdot v}{\longrightarrow} R/L^{[p]}$$ is well-defined.
For instance, if $u\in R$ stands for an element whose image in $R/I^{[p]}$ generates the cyclic module $((I^{[p]}\colon I)\cap (\Omega^{[p]}\colon \Omega))/I^{[p]}$ (see Proposition \ref{useful}) then, clearly, $I$ and $\Omega$ -- as well as $(0)$ and $R$ -- are $u$-compatible ideals.

\begin{Definition}\label{u-clos}\rm Fix some $v\in R$. Given an ideal $L\subset R$, the {\it $\star_v$-closure} of $L$, denoted $L^{\star_v}$,
is the smallest $v$-compatible ideal that contains $L$. When there is no ambiguity as to the element $v$, we simplify $L^{\star_v}$ to $L^{\star}$.

\end{Definition}

\noindent A useful fact about this closure is that $L^{\star}=R$ if and only if $L v  \nsubseteq  {\mathfrak m}^{[p]}$.

\begin{Definition}\label{F-stable}\rm Let $(B, {\mathfrak n})$ be a local ring of dimension $b$ and prime characteristic $p$, 
and let $F\colon H_{\mathfrak n}^b(B)\rightarrow H_{\mathfrak n}^b(B)$ be  the natural Frobenius action. A submodule $N\subset H_{\mathfrak n}^b(B)$ is said to be {\it $F$-stable} if $F(N)\subset N$.

\end{Definition}

Next we recall a couple of key ingredients that we will employ in the proof of Theorem \ref{main}.

\begin{Lemma}\label{Thm1.1}{\rm (De Stefani--N\'u\~nez-Betancourt \cite{DeS-NB})} Let $B$ be an excellent local ring of positive characteristic and satisfying $S_2$ $($e.g., if $B$ is Cohen-Macaulay$).$ Suppose that $B$ has a canonical ideal $K$ such that $B/K$ is $F$-rational. Then, $B$ is strongly $F$-regular.

\end{Lemma}

\begin{Lemma}\label{Smith}{\rm (Smith \cite{S})} Let $(B, {\mathfrak n})$ be an excellent Cohen-Macaulay local ring of dimension $b$ and positive characteristic. Then, the ring $B$ is $F$-rational if and only if  $H_{\mathfrak n}^b(B)$ has no proper non-trivial $F$-stable submodule.

\end{Lemma}

The theorem below gives a sufficient condition for the strong $F$-regularity of our $d$-dimensional ring $A=R/I$ with canonical ideal $\Omega/I$ (we maintain the setup and notation as in the beggining of the section).

\begin{Theorem}\label{main} Suppose that $\Omega \neq R$ and $d\geq 2$. Let $J\subset R$ define the singular locus of $R/\Omega$.
If $$J \, (\Omega^{[p]}\colon \Omega) \, \nsubseteq \, {\mathfrak m}^{[p]}$$ then $A$ is strongly $F$-regular.

\end{Theorem}
\proof Define $B=A/\omega \simeq R/\Omega$. By Lemma \ref{can}, the ring $B$ is Gorenstein of dimension $b=d-1\geq 1$.
According to Proposition \ref{useful}, there is an element $v\in R$ whose image in $R/\Omega^{[p]}$ generates $(\Omega^{[p]}\colon \Omega)/\Omega^{[p]}$ as a cyclic module and such that
$v T$ is the natural Frobenius map on $H^{b}_{{\mathfrak m}}(B)$.

Assume that $$J \, (\Omega^{[p]}\colon \Omega) \, \nsubseteq \, {\mathfrak m}^{[p]}.$$
By Lemma \ref{Thm1.1}, it suffices to prove that the (non-zero) ring $B$ is $F$-rational.
We can identify
$$H^{b}_{{\mathfrak m}}(B) \, = \, \Ann_{E}(\Omega)$$
so that a proper non-trivial $F$-stable submodule $N\subset H^{b}_{{\mathfrak m}}(B)$ can be written as
$N = \Ann_{E}(K)$ for some $v$-compatible ideal $K$ such that $\Omega \subsetneq K\subsetneq R$. What we shall prove is that there is no such $K$ (hence no such $N$). This, together with Lemma \ref{Smith}, will imply the $F$-rationality of $B$.

We will prove, more precisely, that the smallest $v$-compatible ideal that strictly contains $\Omega$ is $R$. To see this, the first step is to invoke the fact described by means of a general algorithm in \cite{K-S}, which guarantees that the smallest {\it radical} $v$-compatible ideal which properly contains $\Omega$ is equal to $\sqrt{L^{\star}}$ ($\star$-closure taken with respect to $\Omega$), where\, $L:= J  \cap  J^\prime$ and
$$J^\prime \, =\, \Ann_R\, \left(\frac{\Omega^{[p]}\colon \Omega}{\Omega^{[p]} + v R}\right).$$
But, by the very nature of $v$, we have $\Omega^{[p]}\colon \Omega=\Omega^{[p]}+ v R$,
which yields $J^\prime=R$ and hence $L=J$.
Note that the hypothesis $J(\Omega^{[p]}\colon \Omega) \nsubseteq  {\mathfrak m}^{[p]}$ is  equivalent to
$$J v \, \nsubseteq \, {\mathfrak m}^{[p]}$$
which in turn means that $J^{\star}=R$, i.e., $\sqrt{L^{\star}}=R$. We will show that every $v$-compatible ideal is radical, and this will finish the proof.

The condition $Jv  \nsubseteq  {\mathfrak m}^{[p]}$ forces $v\notin {\mathfrak m}^{[p]}$, which implies the injectivity of the Frobenius map $v T\colon E\rightarrow E$. To see this, set
$$\Kc \, =\, \{\sigma \in E~|~(v T)(\sigma)=0\}$$
which is the largest submodule of $E$ killed by $v T$.
Moreover, the inclusion $\Kc \subset E$ is compatible with the Frobenius action.
It follows that $$\Kc \, = \, \Ann_{E}(I_1(v))$$ where $I_1(v)$ is the smallest ideal whose Frobenius power contains  $v$.
But $v\notin {\mathfrak m}^{[p]}$, so that $I(v)=R$ and therefore $\Kc=0$, as needed.

Finally, by \cite[Corollary 3.7]{Sh} (see also  \cite[Theorem 3.6]{E-H}), the injectivity of $v T$ implies that all $v$-compatible ideals are radical.
\qed

\begin{Corollary}\label{punctured} Assume that $\Omega \neq R=k[\![x_1, \dots, x_n ]\!]$, and that $R/\Omega$ is locally regular on the punctured spectrum $($i.e., $R_{\fp}/\Omega_{\fp}$ is regular for every $\fp \neq \fm$$)$. If $d\geq 2$ and $$\Omega^{[p]}\colon \Omega \, \nsubseteq \, ({\mathfrak m}^{[p]},\, (x_1\cdots x_n)^{p-1})$$ then $A$ is strongly $F$-regular.

\end{Corollary}
\proof In this situation we have $J={\fm}$, and, moreover, as $x_1, \dots, x_n$ is a regular sequence, $({\mathfrak m}^{[p]}, (x_1\cdots x_n)^{p-1})={\mathfrak m}^{[p]}\colon \fm$. Now the result follows from Theorem \ref{main}. \qed

\begin{Remark}\label{pure}\rm Fedder's criterion of $F$-purity (cf. \cite[Theorem 1.12]{Fe}) states that the ring $A=R/I$ is $F$-pure if and only if\, $I^{[p]}\colon I  \nsubseteq  {\mathfrak m}^{[p]}$. Since strong $F$-regularity is known to imply $F$-purity
(see Section \ref{Section: Preliminaries}), Theorem \ref{main} yields in particular that $A$ is $F$-pure.
Moreover, we point out the similarity between our result and the following characterization given by Glassbrenner in \cite[Theorem 2.3]{Glass}: $A$ is strongly $F$-regular if and only if $$g \, (I^{[p^e]}\colon I) \, \nsubseteq \, {\mathfrak m}^{[p^e]}$$ for {\it every} $A$-regular element $g\in R$ and all $e\gg 0$. Notice that our Theorem \ref{main} requires the existence of {\it some} $g\in J$ satisfying $gv \notin {\mathfrak m}^{[p]}$ (with $v$ as in the proof), while on the other hand it only provides sufficiency -- not necessity in general, as we will illustrate in Example \ref{false-converse}.

\end{Remark}

\begin{Example}\label{symm} \rm Set $A=R/I=k[[x, y, z, w, s, t]]/I_2(\phi)$, where $k$ has characteristic $p=2$ or $p=3$, and $\phi$ is the generic symmetric matrix
$$\phi ~ = ~\left(\begin{array}{ccc}
x  & w & t\\
w  & y & s\\
t  & s & z
\end{array}\right).$$
The ring $A$ is a $3$-dimensional Cohen-Macaulay normal domain with canonical ideal $\Omega/I=(Q+I)/I$, where $Q$ can be taken as the ideal generated by the variables from the first row of $\phi$ (cf.~ \cite{Go}), so that
$$\Omega \, = \, Q\, +\, I\, =\, (x,\, w,\, t,\, yz-s^2).$$ It is easy to see that the ring $R/\Omega$ is locally regular on the punctured spectrum. 
Furthermore, as in this case $\Omega$ is a complete intersection, we can express $v$ simply as the $(p-1)$th power of the product of the generators of $\Omega$. Explicitly,

\[ v\,=\,\left\{ \begin{array}{cr}
    \hspace{-1.5in} xyzwt \, + \,  xwts^2, ~~ p=2\\\\
    \hspace{0.05in} x^2y^2z^2w^2t^2 \, + \, x^2yzw^2t^2s^2 \, + \, x^2w^2t^2s^4, ~~  p=3.

      \end{array} \right.
  \]

\noindent Since $xyzwt$ (resp. $x^2y^2z^2w^2t^2$) lies outside the monomial ideal ${\mathfrak m}^{[2]}$ (resp. ${\mathfrak m}^{[3]}$), we get

\[ v\,\notin\,\left\{ \begin{array}{cr}
    \hspace{-0.19in} ({\mathfrak m}^{[2]}, xyzwts), ~~ p=2\\\\
    \hspace{0.05in} ({\mathfrak m}^{[3]}, (xyzwts)^2), ~~  p=3.

      \end{array} \right.
  \]
By Corollary \ref{punctured}, $A$ is strongly $F$-regular.
\end{Example}

Next, we illustrate that our main hypothesis, $J  (\Omega^{[p]}\colon \Omega) \nsubseteq  {\mathfrak m}^{[p]}$, cannot be relaxed; in particular, it cannot be weakened to the condition $\Omega^{[p]}\colon \Omega \nsubseteq  {\mathfrak m}^{[p]}$.

\begin{Example}\label{N2}\rm Consider the $2$-dimensional Cohen-Macaulay reduced ring $A=R/I=k[[x, y, z, w, s]]/I_2(\psi)$, with $k$ of characteristic $2$, and
$$\psi ~ = ~\left(\begin{array}{cccc}
x  & y & y & s\\
w  & w & z & x
\end{array}\right).$$
A canonical ideal for $A$ is $\Omega/I$ where $\Omega = (x, w, s, yz)$. Therefore $v  =  xyzws$, and
$$J \, (\Omega^{[2]}\colon \Omega) \, = \, \fm \, (\Omega^{[2]},\, v) \, \subset \,  {\mathfrak m}^{[2]}.$$
According to  \cite[Section 9]{K}, the ring $A$ is not $F$-rational, hence it cannot be strongly $F$-regular. Notice, however, that  since $v \notin  {\mathfrak m}^{[2]}$ we have $$\Omega^{[2]}\colon \Omega \, \nsubseteq \, {\mathfrak m}^{[2]}.$$
\end{Example}

Finally, we show that the converse of Theorem \ref{main} is not true.

\begin{Example}\label{false-converse}\rm Let $B=S/L=k[x, y, z, w]/I_2(\varphi)$ where $k$ has characteristic $p=3$ and $\varphi$ is the matrix
$$\varphi ~ = ~\left(\begin{array}{ccc}
x^2  & y & w\\
z  & x^2 & y-w
\end{array}\right).$$
The ring $B$ is a $2$-dimensional (non-Gorenstein) Cohen-Macaulay normal domain, which, according to  \cite[Proposition 4.3]{Singh}, is $F$-regular in the sense that all of its localizations are weakly $F$-regular.

Now notice that $B$ is graded with the grading inherited from the polynomial ring $S$ given by ${\rm deg}(x)=1$, ${\rm deg}(y)={\rm deg}(z)={\rm deg}(w)=2$. But $F$-regular $F$-finite positively graded rings are known to be strongly $F$-regular (\cite[Theorem 2.2(5)]{Singh}). It follows that $B$ is strongly $F$-regular and hence so is its completion $A  = \widehat{B}  =  R/LR =  k[[x, y, z, w]]/I_2(\varphi)$ (see, e.g., \cite[Proof of Theorem 4.1, page 3163]{LS}). We claim, however, that $$\Omega^{[3]}\colon \Omega \, \subseteq \, {\mathfrak m}^{[3]}$$ which will readily illustrate that the converse of our result fails. Indeed, in this case a canonical ideal of $A$ is $\Omega /LR$, where $\Omega = (x^4, y, w)\subset R$ (a complete intersection). Hence $\Omega^{[3]}\colon \Omega=(\Omega^{[3]}, v)$, where $$v \, = \, x^8y^2w^2 \, \in \, {\mathfrak m}^{[3]}.$$

\end{Example}

\section{Generic determinantal rings defined by maximal minors}\label{gen-ring}

We begin this section with a quick recap about determinantal rings defined by (maximal) minors of a matrix of indeterminates over a field. Complete details on the subject can be found in Bruns-Vetter \cite{BV} (cf. also Bruns-Herzog \cite[Section 7.3]{BH}).

We fix a formal power series ring $R=k[[X]]$ over a field $k$, where $X=(x_{ij})_{n\times (n-1)}$ is a matrix of indeterminates, for some $n\geq 2$. Let $$I\,=\,I_{n-1}(X)$$ be the ideal of $R$ generated by the maximal minors of $X$, i.e., its subdeterminants of order $n-1$. It is convenient to write generators explicitly,
$$I\,=\,(G_1, \ldots, G_n)$$ where $G_i$ is the minor which does not involve the $i$th row. By \cite[Theorem 2.1]{BV} the ideal $I$ is a perfect prime ideal of height $2$ 
(indeed, reordering the $G_i$'s and adjusting their signs if necessary, $X$ itself turns out to be a Hilbert-Burch matrix for $I$) and, whenever $i\neq j$, the set $\{G_i, G_j\}$ is a maximal $R$-sequence contained in $I$.

We also consider the ideal $P\subset R$ defined by
\[ 
P=
\left\{ 
\begin{array}{ll}
    R, & \text{if } n=2\\\
    I_{n-2}(X^\prime), & \text{if } n\geq 3\\
      \end{array} \right.
  \]

\noindent where $X^\prime$ is the submatrix consisting of the first $n-2$ rows of $X$ if $n\geq 3$, and in this case $P$ is prime of height $2$ as well. Note that $(G_{n}, G_{n-1})\subset P$, and more precisely, there is a primary decomposition\, $(G_{n}, G_{n-1}) = I \cap  P$.
Furthermore, the canonical ideal of the (normal) Cohen-Macaulay domain $A=R/I$ is $\Omega/I$, where $$\Omega \, = \, P \, + \, I$$ which is also a prime ideal (of height $3$). It is well-known that $A$ is non-Gorenstein if and only if $n\geq 3$ (cf. \cite[Corollary 2.21]{BV}).

\subsection{Explicit generating morphism}

We maintain the preceding setup and notations, and, as in Sections \ref{gen} and \ref{F-rat}, we fix the hypothesis ${\rm char}(k)=p>0$ (but in this part we do not need to assume that $k$ is perfect). Our objective here is to exhibit a generating morphism for the top local cohomology module $$H_I^2(R)~=~H_{I_{n-1}(X)}^2(k[[X]]),$$ 
which is of interest in view of Lyubeznik's theory \cite{L}. Since the ``multiplication by $u$" map
$$\Omega/I \stackrel{\cdot u}{\longrightarrow} \Omega^{[p]}/I^{[p]}$$ is isomorphic to a generating morphism of $H_I^2(R)$ (cf. Remark \ref{generating}), we are reduced to describing $u$ explicitly.

Before revealing the formula for $u$, let us write down the calculations if $n=2$ and $n=3$ (stated below as examples). 
In particular, we will automatically get the base case of the induction used in the proof of Proposition \ref{prod}.

\begin{Example}\label{n=2}\rm The case $n=2$ is easy, and, as mentioned before, is the only situation where $A$ is Gorenstein. We have $R=k[[x, y]]$, $I={\mathfrak m}$ and $\Omega = R$, so that 
$I^{[p]}\colon I=(x^p, y^p, (xy)^{p-1})$ and hence $$u \, = \, (xy)^{p-1} \, = \, (G_2G_1)^{p-1}.$$
\end{Example}

\begin{Example}\label{n=3}\rm We elaborate Fedder's computation in \cite[Proposition 4.7]{Fe} and study the case $n=3$, the first (non-Gorenstein) non-trivial interesting case. Here we have
$$X=\left(\begin{array}{cc}
x  & w\\
y  & s\\
z  & t
\end{array}\right)$$ so that $I=(G_1, G_2, G_3)=(yt-zs,\, xt-zw,\, xs-yw)\subset R=k[[X]]$. Moreover $P=(x, w)$, and hence $\Omega = (G_1, x, w)$.
Since $P$ and $\Omega$ are complete intersections,  $P^{[p]}\colon P = (x^p, w^p, (xw)^{p-1})$ and $\Omega^{[p]}\colon \Omega = (\Omega^{[p]}, (G_1xw)^{p-1})$. Thus
$$\Omega^{[p]}\colon \Omega ~ = ~ (G_1^p,\, x^p,\, w^p,\, (G_1xw)^{p-1}) ~ \subset ~ (G_1^p)\, +\, (P^{[p]}\colon P).$$
Now, pick an arbitrary $f\in (I^{[p]}\colon I)\cap (\Omega^{[p]}\colon \Omega)$. In particular $f\in \Omega^{[p]}\colon \Omega$, and, by the above inclusion, we can write $$f~=~aG_1^p\, +\, e$$ with $a\in R$ and $e\in P^{[p]}\colon P$. For simplicity, set $K=(G_3, G_2)$. Since $K\subset P$, we have $eK\subset P^{[p]}$. Moreover, since $K\subset I$ and $f\in I^{[p]}\colon I$, we get $fK\subset I^{[p]}$. Thus\, $(f - aG_1^p) K  \subset  I^{[p]}$ or, equivalently, $eK\subset I^{[p]}$, and hence $eK\subset I^{[p]}\cap P^{[p]}$. But we know that $K=I\cap P$, which by the exactness of Frobenius gives $$K^{[p]}~=~I^{[p]} \, \cap \, P^{[p]}.$$ 
Therefore $e\in K^{[p]}\colon K=(K^{[p]}, (G_3G_2)^{p-1})$. It follows that there exist $b, c, d\in R$ such that $e=bG_2^p + cG_3^p + d(G_3G_2)^{p-1}$. As $f=aG_1^p + e$, we finally obtain
$$f~ = ~ aG_1^p \, +\, bG_2^p\, +\, cG_3^p\, +\, d (G_3G_2)^{p-1}~ \in ~I^{[p]} \, + \, (u), ~~~~ u\, =\, (G_3G_2)^{p-1}.$$ This proves the inclusion $(I^{[p]}\colon I) \cap  (\Omega^{[p]}\colon \Omega) \subset I^{[p]}+ (u)$.

Now let us verify that, if again we set $u =  (G_3G_2)^{p-1}$, then $I^{[p]}+ (u) \subset (I^{[p]}\colon I) \cap (\Omega^{[p]}\colon \Omega)$, which then will be an equality. 
Since clearly $I^{[p]}\subset \Omega^{[p]}$, we only need to show that $uI\subset I^{[p]}$ and $u\Omega \subset \Omega^{[p]}$. 
Let us prove first that $uI\subset I^{[p]}$. By the exactness of Frobenius we get $${\rm Ass}(R/I^{[p]})\, =\, {\rm Ass}(R/I)$$ which is simply $\{I\}$ since $I$ is prime, 
and hence it suffices to prove that $u=u/1 \in {I_I}^{[p]}\colon I_I$. 
Since $R_I$ is a $2$-dimensional regular local ring, its maximal ideal $I_I$ is a complete intersection of height $2$, and we can write $I_I=(G_3, G_{2})_I\subset R_I$, 
which gives $${I_I}^{[p]}\colon I_I\, =\, ({I}^{[p]},\, (G_3 G_{2})^{p-1})_I\, =\, ({I}^{[p]},\, u)_I$$ and therefore $u \in {I_I}^{[p]}\colon I_I$.

Let us now check that $u\Omega \subset \Omega^{[p]}$. As\, $uI\subset I^{[p]}$ and clearly $\Omega^{[p]}=(P+I)^{[p]}=P^{[p]}+I^{[p]}$, 
it suffices to show that $uP\subset P^{[p]}$, which, since $P$ is prime, amounts to show that $u\in P_P^{[p]}\colon P_P$. We may localize the primary decomposition $(G_{3}, G_{2}) = I \cap  P$ at $P$ and obtain\, $P_P = (G_{3}, G_{2})_P \subset R_P$, hence $P_P^{[p]}\colon P_P=(P^{[p]}, u)_P$, as needed.

Thus we have shown that, in case $n=3$, the cyclic module $((I^{[p]}\colon I)\cap (\Omega^{[p]}\colon \Omega))/I^{[p]}$ is generated by the image of the polynomial $$u \, = \, (G_{3}G_{2})^{p-1}.$$

\end{Example}

\begin{Remark}\rm We point out that, if $n=3$, Fedder \cite[Proposition 4.7]{Fe} explicitly computed the colon ideal $I^{[p]}\colon I$ as being equal to $I^{[p]} + I^{2p-2}$. Here we have not used this fact for our description of the intersection $(I^{[p]}\colon I)\cap (\Omega^{[p]}\colon \Omega)$, and we note that $$(G_{3}G_{2})^{p-1} \, \in \, (I^2)^{p-1}~=~I^{2p-2}.$$ On the other hand, Fedder's computation is no longer valid in higher dimension; for instance if $p=3$ and $n=4$, 
then it can be verified that $$I^{[3]}\colon I \, \neq \, I^{[3]} + I^{4}$$ while Proposition \ref{prod} below computes $(I^{[p]}\colon I)\cap (\Omega^{[p]}\colon \Omega)$ 
for any $n$, so that the behavior illustrated in the examples above is not a coincidence, as we are going to show, in generality, in Proposition \ref{prod}.

\end{Remark}

For convenience, in the present setting we write $$\Cc(X) \, = \, ((I^{[p]}\colon I) \, \cap \, (\Omega^{[p]}\colon \Omega))/I^{[p]}$$ where $X$ is the given $n\times (n-1)$ generic matrix, and as before $u$ denotes a representative for a generator of this cyclic module.

\begin{Proposition}\label{prod} For an arbitrary $n$, we can take\, $u = (G_{n}G_{n-1})^{p-1}$.

\end{Proposition}
\proof We proceed by induction on $n$. 
The case $n=3$ is the content of Example \ref{n=3}.
We set $f=(G_{n}G_{n-1})^{p-1}$ and we will prove that we may take $u=f$ for arbitrary $n$.

Let $Y=(y_{ij})$ be a generic matrix of new indeterminates $y_{ij}$'s over $k$, with $2\leq i \leq n$,\, $2\leq j\leq n-1$, and denote by $S$ the ring obtained by adjoining to the polynomial ring $k[Y]=k[\{y_{ij}\}]$ all the variables that appear in the first row and first column of the original matrix $X$. Now let us invert the variable $x=x_{11}$, that is, we pass to the rings of fractions $S_x$ and $k[X]_x$, where $k[X]=k[\{x_{ij}\}_{n\times (n-1)}]$. It is well-known (cf. \cite[Lemma 7.3.3]{BH}) that the substitution
$$y_{ij} ~ \mapsto ~ x_{ij}\, -\, \frac{x_{1j}\, x_{i1}}{x}, ~~~~ 2\leq i \leq n,~~~ 2\leq j\leq n-1$$ yields a ring isomorphism $\varphi \colon S_x \rightarrow k[X]_x$ such that 
$\varphi(I_{t-1}(Y)_x) = I_t(X)_x$ for every $t\geq 1$. In particular, $$I_x~=~\varphi(I_{n-2}(Y)_x).$$

In analogy with the notation $\Cc(X)$, let $\Cc(Y)$ stand for the corresponding cyclic module derived from the generic matrix $Y$. 
Write generators $$Q \, = \, I_{n-2}(Y) \, = \, (H_2,\ldots, H_n) \, \subset \, S.$$ 
By induction, the class $v+Q^{[p]}$ of the element $$v=(H_nH_{n-1})^{p-1}$$
is a generator of $\Cc(Y)$, 
and  $\varphi(v)$ has the form $f/x^{\alpha}\in k[X]_x$, for some $\alpha\geq 0$. 
Clearly $\varphi(v)\equiv (f/1)~ {\rm mod}\, (I_x)$, and since moreover $\Cc(X)_x\simeq \Cc(Y)_x$ (in virtue of the isomorphism $\varphi$), the residue class $f/1 + I^{[p]}_x$ must generate $\Cc(X)_x$. Therefore, for an arbitrary $g \in k[X]$ such that $g+I^{[p]}$ generates $\Cc(X)$, we get $$(f/1)\equiv (g/1)~ {\rm mod}\, (I^{[p]}_x)$$ or what amounts to the same $$x^{\beta}(f - g)\, \in \, I^{[p]}$$ for some $\beta \geq 0$. Since ${\rm Ass}(R/I^{[p]})=\{I\}$ and $x\notin I$, we necessarily have $f - g\in I^{[p]}$, i.e., the image of $f$ generates $\Cc(X)$ and hence we can take $u=f$, as asserted. \qed

\begin{Remark}\label{u'}\rm For $n\geq 3$, we claim that Proposition \ref{prod} remains valid if we pass from $A=R/I$ to the Gorenstein domain $$A/\omega \, \simeq \, R/\Omega \, =\, R/(I+P)\, =\, k[[X]]/(I_{n-1}(X)+I_{n-2}(X^\prime)).$$ More precisely, if now we denote by $u^\prime \in R$ a representative of a generator of the cyclic $R/\Omega^{[p]}$-module $(\Omega^{[p]}\colon \Omega)/\Omega^{[p]}$, then we can still take $$u^\prime \, = \, (G_{n}G_{n-1})^{p-1}$$ so that the map $R/\Omega \rightarrow R/\Omega^{[p]}$ given by multiplication by this polynomial is isomorphic to a generating morphism of the local cohomology module $H^3_{\Omega}(R)$ (cf. Remark \ref{generating}).

Our claim can be proved by induction in a very similar way, with a suitable adaptation. Let us momentarily assume the validity of the base case $n=3$, and suppose that $n\geq 4$. Let $\varphi \colon S_x\rightarrow k[X]_x$ be the isomorphism described in the proof of Proposition \ref{prod}. Since the submatrix $X^\prime$ (formed by the first $n-2$ rows of $X$) is a generic matrix as well, we may consider the generic submatrix $Y^\prime$ of $Y$ for $i=2,\ldots, n-2$, $j=2,\ldots, n-1$, and we can correspondingly define a ring $S^\prime$ as well as an isomorphism $$\varphi^\prime \colon S^\prime_x\longrightarrow k[X^\prime]_x$$ which is simply the restriction of $\varphi$ to the subring $S^\prime_x\subset S_x$.
Since $I_{n-2}(X^\prime)_x\subset k[X^\prime]_x$ is the image of $I_{n-3}(Y^\prime)_x\subset S^\prime_x$ via $\varphi^\prime$, it follows that the ideal $$P_x \, = \, I_{n-2}(X^\prime)_x \, \subset \, k[X]_x$$ is the extension of the ideal $\varphi^\prime(I_{n-3}(Y^\prime)_x)$ to the ring $k[X]_x$, which in turn coincides with $\varphi(I_{n-3}(Y^\prime)S_x)$. Thus we have $$\Omega_x\, =\, I_x\, + \, P_x \, =\, \varphi (I_{n-2}(Y)_x) \, + \, \varphi^\prime(I_{n-3}(Y^\prime)_x)k[X]_x \, =\, \varphi(I_{n-2}(Y)_x \, +\, I_{n-3}(Y^\prime)S_x).$$ Writing generators $$Q^\prime \, = \, I_{n-3}(Y^\prime)\, +\, I_{n-2}(Y) \, = \, (h_2,\ldots, h_{n-1}, H_2, \ldots, H_n) \, \subset \, S$$ we get, by induction, that if we set
$$v^\prime \, = \, (H_nH_{n-1})^{p-1}$$
then the residue class $v^\prime+ {Q^\prime}^{[p]}$ generates the cyclic module $({Q^\prime}^{[p]}\colon Q^\prime)/{Q^\prime}^{[p]}$, and by construction it satisfies $\varphi (v^\prime)=(G_{n}G_{n-1})^{p-1}/x^{\gamma}$ for some $\gamma \geq 0$. Since $$\Omega_x \, \simeq \, Q^\prime_x$$ via $\varphi$, and since $\Omega$ is prime, the rest of the proof is exactly as in the proof of the proposition.

It remains to check the case $n=3$. We use the same notation as in Example \ref{n=3}, where we verified that
$\Omega^{[p]}\colon \Omega = (\Omega^{[p]}, (G_1xw)^{p-1})$ and $(I^{[p]}\colon I) \cap  (\Omega^{[p]}\colon \Omega) = (I^{[p]}, u)$ for $u=(G_3G_2)^{p-1}$. Denote $u^\prime=(G_1xw)^{p-1}$, so that $\Omega^{[p]}\colon \Omega = (\Omega^{[p]}, u^\prime)$. Since $u\in (I^{[p]}, u) \subset \Omega^{[p]}\colon \Omega$, we get an inclusion $(\Omega^{[p]}, u)\subset (\Omega^{[p]}, u^\prime)$ which turns out to be necessarily an equality since the elements $u, u^\prime$, regarded as homogeneous polynomials of the standard graded ring $k[X]$, have the same degree. Therefore, the class $u^\prime+\Omega^{[p]}=u+\Omega^{[p]}$ generates $(\Omega^{[p]}\colon \Omega)/\Omega^{[p]}$.

\end{Remark}

\subsection{Strong $F$-regularity}

As a consequence of Theorem \ref{main} and Remark \ref{u'}, we will establish a simple proof of the fact that the ring defined by the maximal minors of an $n\times (n-1)$ matrix $X=(x_{ij})$ of 
indeterminates over a perfect field $k$ of characteristic $p\geq 5$, for any $n\geq 2$, is strongly $F$-regular. As mentioned in the Introduction, this is well-known and has been proven in more generality, but the method we have employed in this paper is completely different from the ones available in the literature and we expect that some variation of it may shed a new light on further developments, mainly with a view to the possibility of investigating other classes of determinantal rings.

\begin{Corollary}\label{generic} $($$p\geq 5$$)$ The generic determinantal ring $A=R/I=k[[X]]/I_{n-1}(X)$ is strongly $F$-regular.

\end{Corollary}
\proof This statement is obvious if $n=2$ as in this case we simply have $A\simeq k$. 
Thus we may assume that $n\geq 3$, hence $\dim A = n(n-1)-2\geq 2$.
Let $\Omega$ and $J$ be the ideals as in the statement of Theorem \ref{main} and let $u^\prime$ be as in Remark \ref{u'}.
In order to prove the result by means of Theorem \ref{main}, we need to find some $\Delta \in J$ such that $$\Delta u^\prime \, \notin \, {\mathfrak m}^{[p]}$$ where 
the image of $u^\prime$ in $R/\Omega^{[p]}$ generates $(\Omega^{[p]}\colon \Omega)/\Omega^{[p]}$ as a cyclic module, and $\Omega = I+P=I_{n-1}(X)+I_{n-2}(X^\prime)$.

We treat first the case $n=3$, and we follow the same notation of Example \ref{n=3}. We have $\Omega = (G_1, x, w)=(yt-zs, x, w)$, so that $J=\fm$. According to the fact proved in Remark \ref{u'},
we can take $u^\prime=(G_3G_{2})^{p-1}$. Since $$G_3G_2 \, = \, (xs-yw)(xt-zw) \, = \, -xzws \, + \, \ldots$$ we get, by expansion,
$$u^\prime~=~(-1)^{p-1}(xzws)^{p-1}~+~\ldots$$ Taking for example $\Delta = y$, we obtain that $\Delta u^\prime  \notin  {\mathfrak m}^{[p]}$ and hence $A$ is strongly $F$-regular in this case.

Now suppose that $n\geq 4$. Let $g_j$ denote the determinant of $X^\prime$ with its $j$th column deleted, for $j=1,\ldots, n-1$. 

Then $P=(g_1, \ldots, g_{n-1})$, and note that, for each $j$,
$$\frac{\partial g_j}{\partial x_{ij}}~=~0,~~~ \forall \, i \, \in \, \{1,\ldots, n-2\}.$$
\noindent Since $G_n, G_{n-1}\in P$, we have (minimal) generators $$\Omega \, = \, (G_1, \ldots, G_{n-2}, g_1, \ldots, g_{n-1})$$ whose Jacobian matrix we denote by $\Jc$. We will find an appropriate element $\Delta$ in the ideal $I_3(\Jc)\subset J$ generated by the minors of order $3$ of $\Jc$ (recall that $\Omega$ has height $3$). First, we define the auxiliary monomial
\[ M\,=\,
\left\{ \begin{array}{ll}
     1, & n=4\\
     x_{34}\cdots x_{n-2, n-1}, & n\geq 5.
      \end{array} \right.
  \]
\noindent It is easy to see that
\[ \left\{ \begin{array}{cr}
    \hspace{0.0in} g_1\,=\,x_{12}x_{23}M\,+\, \ldots \\\
    \hspace{-0.06in} g_2\,=\,x_{11}x_{23}M\,+\, \ldots

      \end{array} \right.
  \]
\noindent and expansion along the first column of the $(n-2)\times (n-2)$ matrix whose determinant is $g_3$ yields
$$g_3~=~-x_{21}(x_{12}x_{34}\cdots x_{n-2, n-1}\,+\, \ldots) \, +\, \ldots ~ =~ -x_{21}x_{12}M \, +\, \ldots$$
Now we take the minor $\Delta \in I_3(\Jc)$ given by
$$\Delta~=~{\rm det}\left(\begin{array}{ccc}
\partial g_1/\partial x_{11}  & \partial g_1/\partial x_{21} & \partial g_1/\partial x_{23}\\
\partial g_2/\partial x_{11}  & \partial g_2/\partial x_{21} & \partial g_2/\partial x_{23}\\
\partial g_3/\partial x_{11}  & \partial g_3/\partial x_{21} & \partial g_3/\partial x_{23}
\end{array}\right)~=~{\rm det}\left(\begin{array}{ccc}
0  & 0 & \partial g_1/\partial x_{23}\\
\partial g_2/\partial x_{11}  & \partial g_2/\partial x_{21} & \partial g_2/\partial x_{23}\\
\partial g_3/\partial x_{11}  & \partial g_3/\partial x_{21} & 0
\end{array}\right)$$ that is, $$\Delta~=~\frac{\partial g_1}{\partial x_{23}}\frac{\partial g_2}{\partial x_{11}}\frac{\partial g_3}{\partial x_{21}}~-~\frac{\partial g_1}{\partial x_{23}}\frac{\partial g_2}{\partial x_{21}}\frac{\partial g_3}{\partial x_{11}}.$$ But
$$\frac{\partial g_1}{\partial x_{23}}\frac{\partial g_2}{\partial x_{11}}\frac{\partial g_3}{\partial x_{21}}~=~(x_{12}M)(x_{23}M)(-x_{12}M)~+~\ldots$$
and consequently $$\Delta ~=~-x_{12}^2x_{23}M^3\, +\, \ldots$$

Finally, we claim that $\Delta u^\prime  \notin  {\mathfrak m}^{[p]}$. To prove this, recall that we can take $u^\prime=(G_nG_{n-1})^{p-1}$ (cf. Remark \ref{u'}). Denote by $X_i$ the square submatrix of $X$ obtained by deletion of its $i$th row. Looking at $G_n={\rm det}(X_{n})$ by taking expansion along the first row of $X_n$ yields $$G_n~ =~  \ldots \, + \, (-1)^nx_{1, n-1}\delta_{1, n-1}$$ where the minor $\delta_{1, n-1}$, obtained after deletion of the first row and last column of $X_n$, can be written as
$$\delta_{1, n-1}~ =~ x_{21}\cdots x_{n-1,n-2}~ +~ \ldots$$ so that
$$G_n~ =~ (-1)^nx_{1, n-1} (x_{21}\, \cdots x_{n-1,n-2})~ +~ \ldots $$
Also, note that $$G_{n-1}~ = ~{\rm det}(X_{n-1})~ =~ x_{11}\cdots x_{n-2, n-2}x_{n, n-1} ~+~ \ldots$$
Hence $$G_nG_{n-1}~=~ (-1)^n x_{1, n-1} (x_{21}\cdots x_{n-1,n-2}) (x_{11}\cdots x_{n-2, n-2}x_{n, n-1})~ +~ \ldots$$
and therefore $$u^\prime~=~[(-1)^nx_{1, n-1}\,(x_{21} \cdots  x_{n-1,n-2})\, (x_{11}\cdots x_{n-2, n-2}x_{n, n-1})]^{p-1} ~ + ~\ldots$$
Now, considering the product $\Delta u^\prime$ and noticing, by an elementary inspection, that the monomial
$$x_{12}^2x_{23}M^3\,[x_{1, n-1}\,(x_{21} \cdots  x_{n-1,n-2})\, (x_{11}\cdots x_{n-2, n-2}x_{n, n-1})]^{p-1}$$ avoids ${\mathfrak m}^{[p]}$ if the characteristic satisfies $p\geq 5$, we obtain $\Delta u^\prime  \notin  {\mathfrak m}^{[p]}$ for $p\geq 5$, as needed. \qed

\begin{Remark}\rm Together with some of the well-known facts mentioned in Section \ref{Section: Preliminaries}, the corollary above yields that every ideal of $k[[X]]/I_{n-1}(X)$ is tightly closed and, moreover, that this ring is $F$-pure. Furthermore notice that, for $n=3$ (resp. $n=4$), our proof works for every $p\geq 2$ (resp. $p\geq 3$). For arbitrary $n$, a natural question is whether our argument can be adapted (e.g., by choosing $\Delta$ more efficiently) so as to cover also the remaining cases $p=2$ and $p=3$.

\end{Remark}

\subsection{Further remarks}

We conclude the paper with a couple of remarks.

\begin{Remark}\rm It is natural to ask whether the structural result established in Remark \ref{u'} (which was crucial to the proof of Corollary \ref{generic}) can be extended in an analogous manner to the situation of non-maximal minors, i.e., in the case where $$\Omega \,=\, I_t(X)\, +\, I_{t-1}(X^\prime),~~~~ 2\leq t\leq n-2.$$ Unfortunately, the answer is negative at least if $n=4$, $t=2$ and $p=3$, for which we have verified that $u^\prime$ has degree $18$ as an element of the standard graded polynomial ring $k[X]$, and hence it cannot be expressed as $(Q_1\cdots Q_r)^2$ for quadratic homogeneous polynomials $Q_1, \ldots, Q_r\in I$.

On the other hand, regardless of the ability to understand the entire shape of $u^\prime$, we claim that in this case the ring $A=R/I$ is strongly $F$-regular. In fact, a computation shows that
$$u^\prime\, =\, (x_{11}x_{12}x_{13}x_{22}x_{23}x_{31} x_{33}x_{41}x_{42})^2   \, +\, \ldots $$
and moreover that $x_{43}\in J$, and hence clearly $$x_{43}u^\prime \, \notin \, {\mathfrak m}^{[3]}$$ which, by Theorem \ref{main}, proves the claim.

\end{Remark}

\begin{Remark}\rm We can also raise the problem as to whether our method for finding $u^\prime$, given in Remark \ref{u'}, extends to other classes of determinantal singularities, especially the one formed by generic {\it symmetric} determinantal rings. We have checked, however, that the formula for $u^\prime$ as the $(p-1)$th power of a suitable product of minors does not hold in this situation. For instance, if $I$ is the ideal of minors of order $2$ of a $3\times 3$ symmetric matrix of indeterminates over a field of characteristic $p=2$, then, as we have seen in Example \ref{symm}, the element $u^\prime$ has degree $5$ and consequently it cannot be expressed as a product of quadratic polynomials, and, analogously, if $p=3$ then $u^\prime$ has degree $10$ and hence it cannot be written as the second power of such a product. Of course, even for arbitrary $p$, the shape of $u^\prime$ in the $3\times 3$ case can be easily detected as $\Omega$ is a complete intersection in this situation. However, the problem becomes rather subtle in the case of an $n\times n$ generic symmetric matrix with $n\geq 4$, as the structure of an adequate $u^\prime$ remains quite mysterious and we do not even know whether it can be taken reducible.

Furthermore, inspired by Corollary \ref{generic}, it seems relevant to complement the problem above with the following question:

\begin{Question}\rm Are generic symmetric determinantal rings strongly $F$-regular?
\end{Question}

As far as we know, this is an open problem. The answer is affirmative for $3\times 3$ generic symmetric matrices and $p=2$ or $p=3$, as shown in Example \ref{symm}. There is computational evidence that this is also true in the $4\times 4$ case, at least in low characteristic as well.

Naturally, besides the symmetric case, it would be also of interest to investigate $u^\prime$ as well as the strong $F$-regularity property for other important classes of rings, such as those defined by Pfaffians of generic alternating matrices.

\end{Remark}

\end{document}